\documentclass{siamart}

\usepackage[margin=1in]{geometry}
\usepackage{caption}
\usepackage{multirow}
\usepackage{subfig}

\usepackage{lipsum}
\usepackage{amsfonts}
\usepackage{graphicx}
\usepackage{epstopdf}
\usepackage{algorithmic}
\ifpdf
  \DeclareGraphicsExtensions{.eps,.pdf,.png,.jpg}
\else
  \DeclareGraphicsExtensions{.eps}
\fi

\renewcommand{\vec}[1]{\mathbf{#1}}

\newcommand{\TheTitle}{Algebraic Multigrid Preconditioners for Multiphase Flow in Porous Media} 
\newcommand{\TheAuthors}{Q. Bui, H. Elman, and J.D. Moulton}

\headers{\TheTitle}{\TheAuthors}

\title{{\TheTitle}}%\thanks{This work was funded by the Fog Research Institute under contract no.~FRI-454.}}

\author{
  Quan Bui\thanks{Corresponding author. Applied Math, Stats, and Scientific Computation, University of Maryland, College Park, MD (\email{mquanbui@math.umd.edu}).}
  \and
  Howard Elman\thanks{Department of Computer Science, University of Maryland, College Park, MD (\email{elman@cs.umd.edu}).}
  \and
  David Moulton\thanks{Applied Math and Plasma Physics Group, T-5, Los Alamos National Laboratory, Los Alamos, NM (\email{moulton@lanl.gov}).}
}

\usepackage{amsopn}
\DeclareMathOperator{\diag}{diag}

\ifpdf
\hypersetup{
  pdftitle={\TheTitle},
  pdfauthor={\TheAuthors}
}
\fi

\begin{document}

\maketitle

% REQUIRED
\begin{abstract}
Multiphase flow is a critical process in a wide range of applications, 
including carbon sequestration, contaminant remediation, and groundwater 
management. Typically, this process is modeled by a nonlinear system of partial differential equations derived by considering the mass conservation of each phase (e.g., oil, water), along with constitutive laws for the relationship of phase velocity to phase pressure. In this study, we develop and study efficient solution algorithms for solving the algebraic systems of equations derived from a fully coupled and time-implicit treatment of models of multiphase flow. We explore the performance of several preconditioners based on algebraic multigrid (AMG) for solving the linearized problem, including ``black-box'' AMG applied directly to the system, a new version of constrained pressure residual multigrid (CPR-AMG) preconditioning, and a new preconditioner derived using an approximate Schur complement arising from the block factorization of the Jacobian. We show that the new methods are the most robust with respect to problem character as determined by varying effects of capillary pressures, and we show that the block factorization preconditioner is both efficient and scales optimally with problem size.
\end{abstract}

\section{Introduction}
Multiphase flow is a feature of many physical systems and models of it are used in applications such as reservoir simulation, carbon sequestration, ground water management and contaminant transport. Modeling multiphase flow in highly heterogeneous media with complex geometries is difficult, especially when realistic processes such as capillary pressure are included. The system describing multiphase flow consists of nonlinear partial differential equations, constitutive laws and constraints. In this paper, we focus on the iterative solution of linear systems arising in a fully implicit cell-centered finite volume discretization of single component isothermal two-phase flow model with capillary pressure. This fully implicit time-stepping scheme is among the most robust for simulation of subsurface flow. Moreover, it can serve as a basis for modeling more complex processes in which the physical quantities are tightly coupled. This additional complexity could include adding more components, miscibility between components, thermal effects, and phase transitions.

\indent The fully implicit discretization gives rise to a nonlinear system of equations at each time step. We employ a variant of Newton's method with an exact Jacobian of the discretized equations to solve this system. For the linear system, we use a preconditioned generalized minimal residual (GMRES) method. There is a vast literature on different approaches to precondition the Jacobian system. A very popular approach is to use incomplete LU factorization (ILU) for constructing the preconditioner. Though popular for its general applicability, ILU-based preconditioners are neither effective nor scalable in many cases. Another approach is to consider decoupled preconditioners for the coupled system 
%by Behie and Vinsome \cite{Behie82}. 
\cite{Behie82}. 
This methodology is based on a direct solution of the decoupled pressure system, followed by an iterative solution using ILU for the global system. 
%Wallis \cite{Wallis85} later refined the formulation and proposed solving 
This formulation was refined in \cite{Wallis85}, where it was proposed to solve
the pressure system iteratively, giving rise to the decoupled implicit pressure explicit saturation (IMPES) preconditioner. The effect of the decoupling is to weaken the coupling between pressure and saturation. Thus, it is often used as a preprocessing step to produce a modified Jacobian system, for which new preconditioners can be developed \cite{Clees10,Stueben07}. With recent development of algebraic multigrid (AMG) algorithms, the pressure block can be solved efficiently using AMG, resulting in the constrained pressure residual multigrid (CPR-AMG) approach. In recent developments, AMG has also been applied to solve the coupled system with some success \cite{Clees10,Stueben07}, although developing a general AMG algorithm for these types of problems remains a topic of ongoing research \cite{Osei-Kuffuor15}. Since the Jacobian matrix has a block structure, one can also consider a block LU decomposition with an approximate Schur complement, which has been successfully applied to other models of fluid dynamics \cite{Patankar72,White11}. Besides AMG-based methods, geometric multigrid has also been applied successfully to solve these types of problems \cite{Bastian99,BastianHelmig99}. Our focus in this study is on methods based on AMG because of its general applicability. \\
\indent In this paper, we develop a new block preconditioner designed to respect the coupling inherent in models of multiphase flow, and we report our experience with the performance and scalability of four different preconditioning strategies: (1) a direct algebraic multigrid (AMG) preconditioner for the global system, (2) two-stage CPR-AMG method with correction for the pressure block, also known as the combinative two-stage approach, (3) CPR-AMG with corrections for both the pressure and saturation blocks, known as the two-stage additive approach, and (4) the block factorization (BF) preconditioner. An outline of the paper is as follows. In section 2, we present the mathematical formulation for two-phase flow in porous media and discretization schemes. In section 3, we describe the solution algorithms for the linearized system.
%The
Numerical results for the algorithms are presented in section 4.
We conclude %the this work
with some remarks and discussion of future research directions in section 5.

%\indent The simulator for two-phase flow is implemented within Amanzi, the computational engine of the Advanced Simulation Capability for Environmental Management (ASCEM) project \cite{Amanzi}. Amanzi is interfaced with the Hypre package, a software library for high performance preconditioners and solvers for large, sparse linear systems developed by Lawrence Livermore National Laboratory. In particular, we use BoomerAMG \cite{Henson00} to define our algebraic multigrid preconditioners.

\section{Problem Statement}
We consider isothermal, immiscible two-phase flow through a porous medium. For example, often in reservoir simulation, one phase is oil (the nonwetting phase) and the other is pure water (wetting phase); alternatively, in groundwater management, one may consider a system of contaminated water that infiltrates a domain saturated with air. \\
\indent Conservation of mass of each of the phases leads to the following coupled PDEs:
\begin{align}
& \phi \dfrac{\partial (\rho _w s_w)}{\partial t} + \nabla \cdot (\rho_w \vec{v}_w) = q_w \label{mass_w} \\
& \phi \dfrac{\partial (\rho _n S_n)}{\partial t} + \nabla \cdot (\rho_n \vec{v}_n) = q_n \label{mass_n}
\end{align}
in which $S_w, S_n$ are the saturation, $\rho _w, \rho _n$ are the densities, $q_w, q_n$ are the source terms of the wetting and non-wetting phases respectively, and $\phi$ is the porosity of the medium. We assume a common extension of Darcy's law to multiphase flow and express the phase velocities $\vec{v}_w, \vec{v}_n$ as
\begin{align}
\vec{v}_\alpha = - \dfrac{k_{r\alpha}\vec{K}}{\mu _\alpha}(\nabla P_\alpha - \rho _\alpha g\nabla D), \hspace{5mm} \alpha = w,n. \label{darcys_law}
\end{align}
Here, $\vec{K}$ is the absolute permeability tensor. The terms $k_{r\alpha}, \mu _\alpha, P_\alpha$ are the relative permeability, viscosity, and pressure of phase $\alpha$ respectively, $g$ is the gravitational constant, and $D$ is the depth. We also define the phase mobility $\lambda _\alpha = k_{r\alpha}/\mu _\alpha$. To close the system, we also have the following constitutive law and constraint
\begin{align}
&Pc(S_w) = P_n - P_w \label{capillary_pressure} \\
&S_w + S_n = 1 \label{constraint}
\end{align}
\indent From equations \cref{mass_w,mass_n}, one can derive separate equations for pressure and saturation. The pressure equation is elliptic or parabolic (diffusion); the saturation equation is hyperbolic or convection-dominated. The pressure equation is solved implicitly, and depending on the time discretization strategies applied to the saturation equation, several methods have been developed. In the case where the saturation equation is discretized using an explicit method (e.g., forward Euler), it is referred as IMPES (implicit pressure explicit saturation)\cite{Aziz79}; for an implicit time discretization of the saturation equation, the method is known as the sequential approach, which was first applied to the black-oil model by Watts in 1985 \cite{Watts85}. \\
\indent The appeal of these methods lies in the sequential decoupling between pressure and saturation variables. Each equation can be solved separately. In addition, knowing the characteristics of each equation facilitates the design of efficient preconditioners, which is critical to achieving high performance. Both of these methods have been successfully applied to many problems where the fully implicit method is difficult to implement or shown to be too costly. However, the solution obtained from these approaches may lose accuracy if pressure and saturation are strongly dependent, or if capillary pressure changes very quickly. The lack of accuracy of these methods can be even more pronounced if more complex processes such as miscibility, thermal, and phase transitions are included in the model. For a more complete summary of the advantages and disadvantages of these approaches, we refer to \cite{Lu08}. \\
\indent Substitution of \cref{darcys_law,capillary_pressure} into \cref{mass_w,mass_n} and using the constraint \cref{constraint} yields a system of two equations and two unknowns. Using one popular choice of primary variables, the pressure in the wetting phase and saturation in the nonwetting phase, $\vec{u} = (P_w,S_n)$, we obtain
\begin{align}
-&\dfrac{\partial (\phi\rho _w S_n)}{\partial t} - \nabla \cdot \Big(\rho _w \dfrac{k_{rw}(S_n)}{\mu _w}\vec{K}(\nabla P_w - \rho _w g\nabla D)\Big) = q_w\label{balance1} \\
&\dfrac{\partial (\phi\rho _n S_n)}{\partial t} - \nabla \cdot \Big(\rho _n \dfrac{k_{rn}(S_n)}{\mu _n}\vec{K}(\nabla (P_w + P_c(S_n)) - \rho _n g\nabla D)\Big) = q_n .\label{balance2}
\end{align}
In this paper, we consider solving the coupled system consisting of \cref{balance1,balance2} fully implicitly. We use a cell-centered finite volume method for spatial discretization and the backward Euler method for time discretization, similar to an approach defined in \cite{Dawson97}. This will serve as a base model for adding more complexity in the future. The finite volume method described below is known for its mass conservation property. In addition, it can deal with the case of discontinuous permeability coefficients, and it is relatively straightforward to implement. Under appropriate assumptions, this method also falls into the mixed finite element framework \cite{Peaceman77,Russell83}. For simplicity, we consider a uniform partitioning of the domain $\mathbf{\Omega}$ into 
%uniform 
equal sized
cells $C_i$, 
i.e., $\mathbf{\Omega} = \bigcup\limits_{i=1} C_i$. 
Let $\gamma _{ij}$ denote the area of the face between cells $C_i$ and $C_j$. 
For each cell $C_i$, integration of the mass conservation equations and the divergence theorem gives
\begin{align}
\dfrac{\partial}{\partial t}\int _{C_i} \xi _\alpha + \sum _{j \in \eta _i} \int _{\gamma _{ij}} \psi _\alpha \cdot \vec{n} = \int _{C_i} q_\alpha
\end{align}
where the \textit{storage} $\xi _\alpha = \phi\rho _\alpha S_\alpha$ and the \textit{flux} $\psi _\alpha = \rho _\alpha \vec{v}_\alpha$ terms are approximated using the mid-point rule which is second-order accurate:
\begin{align}
\bar{\xi} _\alpha = \dfrac{1}{V_{C_i}}\int _{C_i} \xi _\alpha, \hspace{3mm}
Q_\alpha = \dfrac{1}{V_{C_i}}\int _{C_i} q_\alpha.
\end{align}
The surface integrals are discretized using two-point flux approximation (TPFA); dropping the phase subscript, this gives
\begin{align}
&\int _{\gamma _{ij}} \psi \cdot \vec{n} = -\gamma _{ij} \Big(\rho \dfrac{k_{r}}{\mu} \vec{K} \Big)_{ij + 1/2} \big( \omega _i - \omega _j \big)\,, \\
&\omega _i = \dfrac{P_i - \rho_{ij+1/2}\; g D_i}{\Delta x_{ij+1/2}}\;.
\end{align}
The index $ij + 1/2$ signifies an appropriate averaging of properties at the interface between cell $i$ and $j$. The coefficients $(\rho k_{r}/\mu)_{ij + 1/2}$ are approximated by upwinding based on the direction of the velocity field, i.e.,
\begin{align}
\Big(\rho \dfrac{k_{r}}{\mu}\Big)_{ij + 1/2} = 
\begin{cases}
\Big(\rho \dfrac{k_{r}}{\mu}\Big)_{i}\;, \hspace{5mm} \text{if } \vec{v} \cdot \vec{n} > 0 \\
\Big(\rho \dfrac{k_{r}}{\mu}\Big)_{j}\;, \hspace{5mm} \text{otherwise}
\end{cases}
\end{align}
and the absolute permeability tensor on the faces are computed using harmonic averaging,
\begin{align}
\vec{K} _{ij} = (\Delta x_i + \Delta x_j)\Big(\dfrac{\mathbb{K}_i\mathbb{K}_j}{\Delta x_i\mathbb{K}_j + \Delta x_j\mathbb{K}_i}\Big).
\end{align}
Discretization in time using the backward Euler method gives a fully discrete system of nonlinear equations,
\begin{align}
(\bar{\xi})_i^{n+1} - (\bar{\xi})_i^n = 
&- \dfrac{\bigtriangleup t}{V_{C_i}}\sum _{j \in \eta _i} \gamma _{ij} \left(\rho \dfrac{k_{r}}{\mu} \mathbb{K}\right)_{ij + 1/2}^{n+1} \Big( \omega _i^{n+1} - \omega _j^{n+1} \Big) - Q^{n+1}. \label{fully_discrete_system}
\end{align}

\section{Solution Algorithms}
The system of nonlinear equations \eqref{fully_discrete_system} can be written generically as $F(u) = 0$ where $F: \mathbb{R}^n \rightarrow \mathbb{R}^n$. We solve the system using Newton's method, which requires solution of a linear system at each iteration $k$:
\begin{align}
\dfrac{\partial F}{\partial u}\Big|_{u=u_k} (u_{k+1} - u_k) = - F(u_k) .
\end{align}
In our case, the solution vector $u$ consists of all the pressure and saturation unknowns at all the cell centers. The Jacobian system resulting from the derivative $\partial F/\partial u$ is often very difficult to solve using iterative methods, and preconditioning is critical for rapid convergence of Krylov subspace methods such as GMRES \cite{Saad86}. Next, we discuss the linear system arising from the Newton's method and give a detailed description of the solution algorithms we will use to solve this system.
\subsection{Linear System}
For the set of primary variables $u = (P_w,S_n)$, assuming unknowns corresponding to physical variables are grouped together and unknowns associated with nodal points in the domain are ordered lexicographically, each nonlinear Newton iteration entails the solution of a discrete version of a block linear system of the form
\begin{align}
\begin{pmatrix}
-\nabla \cdot (\lambda _wK\nabla) & -\dfrac{\phi}{\partial t} - \nabla \cdot (\vec{v}_w )\\
-\nabla \cdot (\lambda _nK\nabla) & \dfrac{\phi}{\partial t} + \nabla \cdot (\vec{v}_n) + \nabla \cdot (\lambda _n P_c^{\prime}K\nabla)
\end{pmatrix}
\begin{pmatrix}
\delta P_w \\
\delta S_n
\end{pmatrix}
= -\begin{pmatrix}
q_w \\
q_n
\end{pmatrix} \label{linearize_fc}
\end{align}
in which 
\begin{align}
&\vec{v}_w = - \lambda _w^{\prime}K\nabla\tilde{P}_w \\
&\vec{v}_n = - \lambda _n^{\prime}K\nabla\tilde{P}_n + \lambda _nK\nabla(P_c^{\prime}) \;.
\end{align}
All the coefficients in equation \eqref{linearize_fc} are evaluated at the
linearization point $\tilde{P}_w,\tilde{S}_n$.
In a more concise form, the Jacobian matrix of the system has $2\times 2$ block structure
\begin{align}
J =
\begin{pmatrix}
A_{pp} & A_{ps} \\
A_{sp} & A_{ss}
\end{pmatrix},
\label{jacobian}
\end{align}
and the linear system is $Jc=q$.
The characteristics of the matrix have been discussed in numerous
papers \cite{BastianHelmig99,Dawson97,Klie96,Stueben07}.
We summarize important characteristics of the operators here:
\begin{itemize}
\item $J$ is nonsymmetric and indefinite
\item The block $A_{pp}$ has the structure of a discrete purely elliptic problem for pressure.
\item The coupling block $A_{ps}$ has the structure of a discrete first-order hyperbolic problem in the non-wetting phase saturation.
\item The coupling block $A_{sp}$ has the structure of a discrete convection-free parabolic problem in the wetting phase pressure.
\item The block $A_{ss}$ has the structure of a discrete parabolic (convection-diffusion) problem for saturation when capillary pressure is a non-constant function of the saturation. When capillary pressure is zero or a constant, there is no diffusion term and the block has the form of a hyperbolic problem.
\item Under mild conditions, i.e. modest time-step size, the blocks $A_{pp},A_{ps},A_{ss}$ are diagonally dominant.
\end{itemize}
In this paper, we present some numerical results that show how different models of capillary pressure affect the algebraic properties of the (2,2)-block $A_{ss}$ in particular and the global system in general, which consequently determines the success of AMG solution algorithms.
Our emphasis is on the development and use of preconditioning operators denoted $M \approx J$,
for the purpose of solving preconditioned systems
\begin{equation} \label{precon-system}
    J M^{-1} \hat c = q, \quad c=M^{-1} \hat c .
\end{equation}

%\subsection{Decoupling Operators}
%Designing effective preconditioners for the Jacobian system \eqref{jacobian} is challenging because of the strong coupling between pressure and saturation. To weaken the coupling between pressure and saturation unknowns, one can often construct a decoupling operator $D$, which acts as a matrix scaling to the original Jacobian
%\begin{align}
%\tilde{J} \equiv D^{-1}J \equiv \begin{pmatrix}
%\tilde{A_{pp}} & \tilde{A_{ps}} \\
%\tilde{A_{sp}} & \tilde{A_{ss}}
%\end{pmatrix}\;.
%\end{align}
%A natural choice for $D$ that promotes code reuse is an IMPES type operator \cite{Dawson97}
%\begin{align}
%D = \begin{pmatrix}
%D_{pp} & D_{ps} \\
%D_{sp} & D_{ss}
%\end{pmatrix}
%= \begin{pmatrix}
%\text{diag}(A_{pp}) & \text{diag}(A_{ps}) \\
%\text{diag}(A_{sp}) & \text{diag}(A_{ss})
%\end{pmatrix}
%\end{align}
%Other approaches include quasi-IMPES \cite{Lacroix03} and alternate block factorization \cite{Bank89} strategies. In general, the decoupling operator applied to the Jacobian system results in clustering of the eigenvalues around one, and this operator can be an effective preconditioner. For a more detailed discussion of the effects of this decoupling step, we refer to \cite{Klie96}. Here however, we do not employ any decoupling strategies, in order to examine the robustness of BoomerAMG for fully coupled systems. Moreover, when more complex processes are added to the model, there is no guarantee that a good decoupling strategy exists.

\subsection{Algebraic Multigrid}
Multigrid is a highly efficient and scalable method available for solving large sparse linear systems \cite{Yang06}. Geometric multigrid uses a hierarchy of nested grids, whose construction depends on the geometry of the problem and \textit{a priori} knowledge of the grids. AMG methods such as those developed in \cite{Stueben01} have the advantage of not requiring an explicit hierarchy of nested grids. AMG constructs coarse grids based on the matrix values only, which makes it suitable for solving a wide range of problems on complicated domains and unstructured grids. Despite its successful application to scalar problems, using AMG for coupled systems is still relatively limited. Some attempts to use AMG to solve fully coupled systems encountered in modeling multiphase flow for reservoir simulation include \cite{Clees10,Stueben07}. In this work, we use BoomerAMG \cite{Henson00}, part of the Hypre package \cite{Falgout06,Falgout02}, as a black-box AMG solver. We note that in order to use BoomerAMG for the coupled system in our case, the Jacobian matrix needs to be ordered by grid points, i.e.
\begin{align}
J = \begin{pmatrix}
A_{11} & \hdots & A_{1N} \\
\vdots & \ddots & \hdots \\
A_{N1} & \hdots & A_{NN}
\end{pmatrix}
\end{align}
in which $N$ is the number of grid points, and $A_{ij}$ are $2\times 2$ matrices representing the couplings between pressure and saturation at points $i$ and $j$. This is called the ``point'' method in \cite{Stueben07}.

\subsection{Two-stage Preconditioning with AMG}
Unlike AMG, which has not been popular in reservoir simulation until recently, two-stage preconditioners are widely used \cite{Klie96}. This idea was developed and first appeared in the context of multiphase flow modeling in the work of Wallis \cite{Wallis85}. Following \cite{Dawson97}, we refer to this method as the constrained pressure residual (CPR) approach. There are many variants of two-stage preconditioners. We discuss two algorithms here: the two-stage combinative preconditioner - CPR-AMG(1), and the two-stage additive preconditioner - CPR-AMG(2) \cite{Axelsson94}. \\
\textbf{Algorithm 1.} Two-stage Combinative - CPR-AMG(1)
\begin{enumerate}
\item At each iteration $k$ let the residual be $r_k$.
\item Solve $u_{k+1/2} = P_1^{-1} r_k$.
\item Update the residual $r_{k+1/2} = r_k - Au_{k+1/2}$.
\item Solve for the pressure correction $A_{pp} \delta _p = R_p r_{k+1/2}$
\item Update the solution $u_{k+1} = u_{k+1/2} + R_p^T\delta _p$
\end{enumerate}
%\vspace{2mm}
\textbf{Algorithm 2.} Two-stage Additive - CPR-AMG(2)
\begin{enumerate}
\item At each iteration $k$ let the residual be $r_k$.
\item Solve $u_{k+1/2} = P_1^{-1} r_k$.
\item Update the residual $r_{k+1/2} = r_k - Au_{k+1/2}$.
\item Solve for the pressure correction $A_{pp} \delta _p = R_p r_{k+1/2}$
\item Solve for the saturation correction $A_{ss} \delta _s = R_s r_{k+1/2}$
\item Update the solution $u_{k+1} = u_{k+1/2} + R_p^T\delta _p + R_s^T\delta _s$
\end{enumerate}
%\vspace{2mm}
The matrices $R_p, R_s$ denote the restriction of the global unknown vector to those associated with pressure and saturation respectively. That is, $R_p \in \mathbb{R}^{n\times 2n}$ and for $u = \begin{pmatrix} p\\ s
\end{pmatrix}$
\begin{align}
R_p u = p, \hspace{5mm} &R_p^{T} u = 
\begin{pmatrix}
p \\
0
\end{pmatrix}; \hspace{11mm} R_s u = s \hspace{5mm} 
R_s^{T} u = 
\begin{pmatrix}
0 \\
s
\end{pmatrix}
\end{align}
Then, in matrix form, the action of the two-stage preconditioners can be expressed as
\begin{align}
&u = M_{comb}^{-1}r = (I - R_p^T A_{pp}^{-1} R_p (A - P_1))P_1^{-1}r \label{m_comb}\\
&u = M_{add}^{-1}r = (I - (R_p^T A_{pp}^{-1} R_p + R_s^T A_{ss}^{-1} R_s)(A - P_1))P_1^{-1}r \label{m_add}
\end{align}
For the preconditioner $P_1$ in step 2 of both algorithms, we use the incomplete factorization with no fill ILU(0) method. For the correction solve, we apply AMG with one V-cycle iteration. The combinative approach with AMG was presented in \cite{Lacroix03,Masson12}. However, this method does not work well in the presence of fast changing capillary pressure. We confirm this observation in the next section. To deal with fast changing capillary pressure, we employ an additive CPR-AMG approach, which involves one more AMG solve for the correction of the saturation block. The intuition is that when the absolute value of the derivative of capillary pressure $|dP_c/dS_w|$ is large, the block $A_{ss}$ becomes diffusion dominated, and AMG can handle %deal with 
it efficiently.

\subsection{Block Factorization Preconditioners}
Consider the following decomposition of the Jacobian,
\begin{align*}
J = \begin{pmatrix}
A_{pp} & A_{ps} \\
A_{sp} & A_{ss}
\end{pmatrix}
= \begin{pmatrix}
I & A_{ps}A_{ss}^{-1} \\
0 & I
\end{pmatrix}
\begin{pmatrix}
S & 0 \\
0 & A_{ss}
\end{pmatrix}
\begin{pmatrix}
I & 0 \\
A_{ss}^{-1}A_{sp} & I
\end{pmatrix},
\end{align*}
where $S$ is the Schur complement 
\begin{align}
S = A_{pp} - A_{ps}A_{ss}^{-1}A_{sp} .
\end{align}
We could choose 
\begin{align}
&M = \begin{pmatrix}
I & A_{ps}A_{ss}^{-1} \\
0 & I
\end{pmatrix}
\begin{pmatrix}
S & 0 \\
0 & A_{ss}
\end{pmatrix}
= \begin{pmatrix}
S & A_{ps} \\
0 & A_{ss}
\end{pmatrix}
\end{align}
as an upper-triangular block preconditioner; 
this incorporates the effects of
%which can capture 
the coupling block $A_{ps}$, which contains the time derivative and gravity terms \eqref{linearize_fc}. We use an approximation of the Schur complement in which $A_{ss}$ is replaced by its diagonal values:
\begin{align}
\tilde{S} = A_{pp} - A_{ps}\diag{(A_{ss})}^{-1}A_{pp}. \label{SIMPLE-Schur}
\end{align}
The purpose of this is to keep the Schur complement sparse so that the action of its inverse can be applied efficiently. This idea is the basis of the SIMPLE method used in other models of fluid dynamics \cite{Patankar72}. A similar approach has also been applied to problems in single phase flow coupled with geomechanics in \cite{White11}.\\
\textbf{Algorithm 3.} Block factorization preconditioner
\begin{enumerate}
\item At each iteration $k$ let the residual be $r_k$.
\item Solve for the saturation $A_{ss} s_{k+1} = R_s r_{k}$ using AMG.
\item Compute the residual for pressure $r = R_p r_k - A_{ps} s_{k+1}$.
\item Solve for the pressure $\tilde{S} p_{k+1} = r$ using AMG.
\end{enumerate}
An important advantage of this algorithm is that it does not rely on 
an ILU factorization.
In matrix form,
\begin{align}
M_{bf}^{-1} = \begin{pmatrix}
\tilde{S}^{-1} & -\tilde{S}^{-1}A_{ps}A_{ss}^{-1}\\
0 & A_{ss}^{-1}
\end{pmatrix}\;. \label{m_bf}
\end{align}

\section{Numerical Results}
In this section, we perform numerical experiments for the four aforementioned preconditioners. All of them are implemented in Amanzi, a parallel open-source multi-physics C++ code developed as a part of the ASCEM project \cite{Amanzi}. Although Amanzi was first designed for simulation of subsurface flow and reactive transport, its modular framework and concept of process kernels \cite{ECoon_JDMoulton_SPainter_2014a} allow new physics to be added relatively easily for other applications. The two-phase flow simulator employed in this work is one such example. Amanzi works on a variety of platforms, from laptops to supercomputers. It also leverages several popular packages for mesh infrastructure and solvers through a unified input file. Here, all of our experiments use a classical AMG solver through BoomerAMG in Hypre. The ILU(0) method is from Euclid, also a part of Hypre. GMRES is provided within Amanzi. For simplicity, we employ structured Cartesian grids for the test cases, but we can also use unstructured K-orthogonal grids. This section has three parts. In the first part, we show the results for a two-dimensional oil-water model problem. Although the problem is small, it is difficult to solve due to the heterogeneity of the permeability field. In the second part, we report the results for a three-dimensional example. In the last part, we examine the scalability of the three preconditioning strategies. Unless specified otherwise, we use the benchmark problem SPE10 \cite{Christie01} for permeability data and porosity.

\subsection{Two-dimensional oil-water problem} \label{2D-oil-water}
\begin{figure}
\centering
\captionsetup{justification=centering}
\caption{Capillary pressure curves for Brooks-Corey model with entry pressure $P_d = 10^5 \hspace{2mm} Pa$.} \label{brooks_corey}
\includegraphics[width=0.5\textwidth]{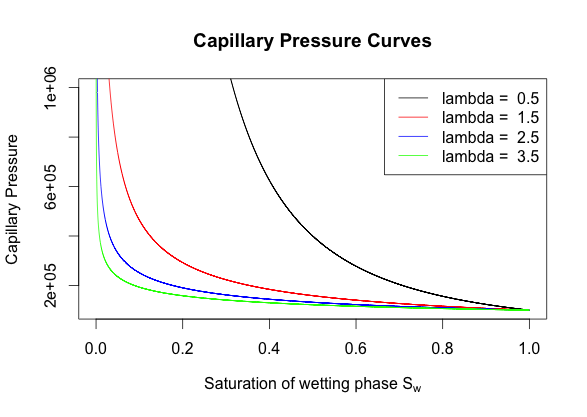}
\end{figure}

The domain is a rectangle with dimensions $762 \times 15.24$ meters. The mesh is $100 \times 20$, which means that the problem is truly two-dimensional in the $xz$ plane. We inject pure water into the domain through the boundary at the lower left corner, and oil and water exit the domain through the top right corner. These correspond to the $S_w = 1.0, \hspace{1mm} \lambda _w \nabla P_w \cdot \vec{n} = -50 \hspace{2mm} m^3/day$ at the inlet, and $S_w = 0.2, \hspace{1mm}P_w = 0$ at the outlet. The simulation is run for 200 days with time step $\Delta t = 20 \text{  days}$. \\
\indent For capillary pressure models, we employ a simple linear model and the Brooks-Corey \cite{BrooksCorey64} model:
\begin{align}
&\text{Linear model: } P_c(S_w) = P_0(1-\bar{S}_{w}), \hspace{5mm}\text{Brooks-Corey: } P_c(S_w) = P_d \bar{S}_{w}^{-1/\lambda} 
%\text{Van Genuchten: } P_c(S_w) = \dfrac{1}{\alpha}\Big( S_{we}^{-1/m} -1 \Big)^{1/n}\.,
\end{align}
in which $\bar{S}_w$ is the effective saturation, $P_d$ is the entry pressure, and $\lambda$ is related to the pore-size distribution. For the Brooks-Corey model, the typical range of $\lambda$ is $\left[ 0.2,3.0\right]$ \cite{Bastian99,Corey94}. In general, $\lambda$ is greater than 2 for narrow distributions of pore sizes, and $\lambda$ is less than 2 for wide distributions. For example, sandpacks with broader distributions of particle sizes have $\lambda$ ranging from $1.8$ to $3.7$ \cite{RelPerm}. The Brooks-Corey capillary pressure curves for various values of $\lambda$ are plotted in Figure \ref{brooks_corey}. Other parameters are listed in Table \ref{params_q5} and example 1 of Table \ref{pc_params}.
\indent For all of the simulations presented here, the convergence tolerance for Newton's method is $||F(x)|| \le 10^{-6}$, and the linear tolerance for GMRES is $||J\delta u_k - F(u_k)|| \le 10^{-12}||F(u_k)||$, which is the default in Amanzi. BoomerAMG is used as a preconditioner. The number of V-cycle steps is set to 1. The coarsening strategy is the parallel Cleary-Luby-Jones-Plassman (CLJP) coarsening \cite{Cleary98}. The interpolation method is the classical interpolation defined in \cite{RugeStueben86}, and the smoother is the forward hybrid Gauss-Seidel / SOR scheme.
\begin{figure}
\centering
\includegraphics[width=0.8\textwidth]{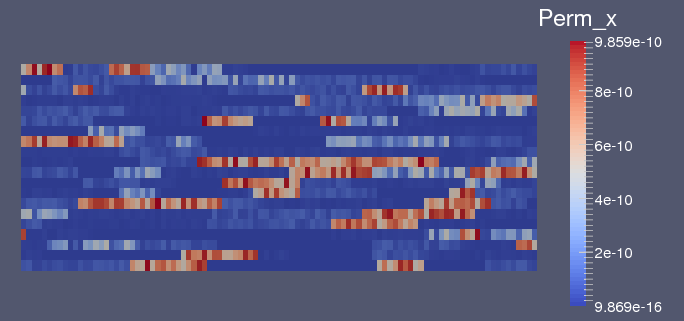}
\captionsetup{justification=centering}
\caption{Permeability field obtained from SPE10 model 1 data. \\
The x-direction is scaled down by $1/20$ for visualization.}
\end{figure}
\begin{table}
%\begin{minipage}[t]{.5\textwidth}
\centering
\caption{Input data for the quarter-five spot problem.} \label{params_q5}
\begin{tabular}{l l}
\hline
Initial wetting phase pressure & $10^5$ $Pa$ \\
Initial nonwetting phase saturation & 0.8 \\
Residual wetting phase saturation & 0.0 \\
Nonwetting phase density & 700 $kg/m^3$\\
Wetting phase density & 1000 $kg/m^3$\\
Nonwetting phase viscosity & 10.0 $cP$\\
Wetting phase viscosity & 1 $cP$\\
Porosity & 0.2 \\
\hline
\end{tabular}
%\end{minipage}%
%\begin{minipage}[t]{0.5\textwidth}
%\centering
%\caption{Parameters for capillary pressure models} \label{pc_params_q5_1}
%\begin{tabular}{l l}
%\hline
%Linear entry pressure $P_0$ & $10^5$ $Pa$ \\
%Brooks-Corey entry pressure $P_d$ & $10^6$ $Pa$ \\
%Brooks-Corey $\lambda$ & 2.5 \\
%\hline
%\end{tabular} \\
%\begin{table}[H]
%\centering
%\caption{Parameters for capillary pressure models} \label{pc_params_q5_2}
%\begin{tabular}{l l}
%\hline
%Linear entry pressure $P_0$ & $10^4$ $Pa$ \\
%Brooks-Corey entry pressure $P_d$ & $10^5$ $Pa$ \\
%Brooks-Corey $\lambda$ & 0.8 \\
%\hline
%\end{tabular}
%\end{table}
%\end{minipage}
\centering
\caption{Parameters for capillary pressure models} \label{pc_params}
\begin{tabular}{|l|c|c|c|c|}
\hline
Parameters & Ex 1 & Ex 2 & Ex 3 & Ex 4 \\
\hline
Linear entry pressure $P_0$ & $10^5$ & $10^4$ & $10^3$ & $10^6$ \\
Brooks-Corey entry pressure $P_d$ & $10^6$ & $10^5$ & $ 2\times 10^4$ & $10^6$ \\
Brooks-Corey $\lambda$ & 2.5 & 0.8 & 2.5 & 0.8 \\
\hline
\end{tabular}
%\vspace{-5mm}
\end{table}
\par In order to explore the effects of different models for capillary pressure on solver performance, we use the four sets of parameters listed in Table \ref{pc_params}. In Example 1, the parameters are chosen such that the $L_\infty$ norm of the derivative of capillary pressure $P_c^\prime$ is large, leading to a diffusion-dominated case (see equation \eqref{linearize_fc}). In Example 2, the parameters are tuned to reduce the $L_\infty$ norm of $P_c^\prime$, leading to an advection-dominated case. Example 3 is a more extreme case of example 2, in which $P_c^\prime$ is further decreased, leading to a strongly advection-dominated case. We also note the difference between the linear model and the Brooks-Corey model for capillary pressure. The derivative $P_c^\prime$ for the linear model is a constant value, which means that the character of the problem, i.e. diffusion-dominated or advection-dominated, is the same everywhere for the whole domain. In the Brooks-Corey model, $P_c^\prime$ depends on the saturation of the wetting phase, and the problem can be diffusion-dominated in one part of the domain, and advection-dominated in another part. This can cause further difficulties for AMG-based solvers, whose optimal performance is sensitive to the characteristics of the problem. 

\begin{table} 
\centering
\captionsetup{justification=centering}
\caption{Performance of three preconditioning strategies for set of parameters in the diffusion-dominated Example 1.} \label{performance_q5_1}
\begin{tabular}{| c | c | c | c | c | c | c | c | c |}
\hline
\multirow{2}{*}{Methods/Models} & \multicolumn{4}{c|}{Linear} & \multicolumn{4}{c|}{Brooks Corey} \\
\cline{2-9}
& NI & LI & LI/NI & Time & NI & LI & LI/NI & Time \\
\hline
AMG & 32 & 368 & 11.5 & 27.2 & 36 & 470 & 13.1 & 37.24\\
CPR-AMG(1) & 32 & 3695 & 115.5 & 324.15 & 36 & 5831 & 162 & 567.7 \\
CPR-AMG(2) & 32 & 899 & 28.1 & 103.94 & 36 & 1102 & 30.6 & 134.6\\
BF & 32 & 524 & 16.4 & 33.17 & 36 & 599 & 16.6 & 46.2\\
\hline
\end{tabular}
%\vspace{-5mm}
\end{table}

The performance of the three strategies is summarized in Tables \ref{performance_q5_1}, \ref{performance_q5_2}, and \ref{performance_q5_3}. NI denotes the number of nonlinear iterations, LI the number of linear iterations, LI/NI the average number of linear iterations per nonlinear iterations, and Time the total time 
in seconds of the whole simulation.   %\\
For the diffusion-dominated problem for which the results are shown in Table \ref{performance_q5_1}, AMG is the most efficient method, about 25\% more efficient than the block preconditioner in terms of both iteration counts (linear iterations per Newton step) and total run time.  Note that for this example, the diffusion terms in the (2,2)-block ($A_{ss}$) are large and the performance of AMG tends to be very strong for problems of this type. For the linear model, the block factorization approach still takes about 8 times fewer linear iterations, and it is about 10 times faster in total run time than CPR-AMG(1). The reason for this discrepancy is that CPR-AMG(1) is a two-stage preconditioner, and it requires an extra global solve using ILU. The block factorization preconditioner does not rely on ILU, which helps improve the run time significantly. CPR-AMG(2) also performs well in this case. Although it requires one more AMG solve per Newton iteration than CPR-AMG(1), it still outperforms CPR-AMG(1) in terms of both the number of linear iterations per Newton step and the total run time. The same conclusion can be made for the Brooks-Corey model.

\begin{table} 
\centering
\captionsetup{justification=centering}
\caption{Performance of three preconditioning strategies for set of parameters in the advection-dominated Example 2.} \label{performance_q5_2}
\begin{tabular}{| c | c | c | c | c | c | c | c | c |}
\hline
\multirow{2}{*}{Methods/Models} & \multicolumn{4}{c|}{Linear} & \multicolumn{4}{c|}{Brooks Corey} \\
\cline{2-9}
& NI & LI & LI/NI & Time & NI & LI & LI/NI & Time\\
\hline
%Time(s) & 149.8 & 225.3 & 238 \\
AMG & 37 & 2575 & 69.6 & 138.8 & - & - & - & -\\
CPR-AMG(1) & 37 & 1919 & 51.9 & 175.5 & 55 & 4851 & 88.2 & 605.7 \\
CPR-AMG(2) & 37 & 1222 & 33.0 & 157.1 & 55 & 3701 & 67.3 & 506.8 \\
BF & 37 & 684 & 18.5 & 51.7 & 55 & 1633 & 29.7 & 131.1\\
\hline
\end{tabular}
\centering
\captionsetup{justification=centering}
\caption{Performance of three preconditioning strategies for set of parameters in the strongly advection-dominated Example 3.} \label{performance_q5_3}
\begin{tabular}{| c | c | c | c | c | c | c | c | c |}
\hline
\multirow{2}{*}{Methods/Models} & \multicolumn{4}{c|}{Linear} & \multicolumn{4}{c|}{Brooks Corey} \\
\cline{2-9}
& NI & LI & LI/NI & Time & NI & LI & LI/NI & Time\\
\hline
AMG & - & - & - & - & - & - & - & -\\
CPR-AMG(1) & 43 & 1079 & 25.1 & 122.8 & 48 & 2173 & 45.3 & 247.6 \\
CPR-AMG(2) & 43 & 1442 & 35.5 & 169.8 & 48 & 4805 & 100.1 & 560.5\\
BF & 43 & 1002 & 23.3 & 69.8 & 48 & 1829 & 38.1 & 121.8\\
\hline
\end{tabular}
\end{table}
The results reported in Table \ref{performance_q5_2} reveal the lack of robustness of AMG when applied to the coupled system. In contrast to the diffusion-dominated case, for the linear model of capillary pressure, AMG requires the highest number of linear iterations per Newton step for the advection-dominated case, and it even diverges for the Brooks-Corey model. The block factorization preconditioner still shows good performance, taking about half the number of iterations and running four times faster than the next best method, which is CPR-AMG(2). CPR-AMG(1) is still the least effective method in this case for both capillary pressure models.

For the strongly advection-dominated problem with parameters in example 3, AMG diverges for both the linear and Brooks-Corey capillary pressure models. The performance of CPR-AMG(2) is also affected in this case, trailing that of CPR-AMG(1). CPR-AMG(2) is still more robust than direct application of AMG, however, since unlike AMG, this method still converges. The block factorization preconditioner is again the most effective method, requiring fewer number of iterations and about half the run time of CPR-AMG(1). This suggests that when the diffusion term in the $A_{ss}$ block gets small, the coupling block $A_{ps}$ which has the structure of a discrete first-order hyperbolic problem for the saturation, becomes important and needs to be taken into account.

\subsection{Two-dimensional problem with gravity}
%The domain is a 2D box with uniform mesh. We increase the size of the domain as we make the mesh larger to keep the refinement $h$ constant.
%\begin{table}[H]
%\centering
%\caption{Iteration counts for diffusion-dominated case.}
%\begin{tabular}{| c | c | c | c | c |}
%\hline
%Methods/Models & $20^2$ & $40^2$ & $80^2$ & $160^2$\\
%\hline
%%Time(s) & 149.8 & 225.3 & 238 \\
%AMG & 16 & 9 & 9 & 9\\
%CPR-AMG(1) & 13.4 & 13.7 & 13.9 & 14.6 \\
%CPR-AMG(2) & 22.2 & 23.5 & 23.7 & 24.1 \\
%BF & 14.1 & 14.1 & 13.7 & 13.3 \\
%\hline
%\end{tabular}
%\centering
%\caption{Iteration counts for advection-dominated case.}
%\begin{tabular}{| c | c | c | c | c |}
%\hline
%Methods/Models & $20^2$ & $40^2$ & $80^2$ & $160^2$\\
%\hline
%%Time(s) & 149.8 & 225.3 & 238 \\
%AMG & - & - & - & -\\
%CPR-AMG(1) & 15.0 & 15.3 & 15.6 & 16.4 \\
%CPR-AMG(2) & 27.9 & 29.8 & 30.4 & 31.7\\
%BF & 9.7 & 9.7 & 9.8 & 10.8\\
%\hline
%\end{tabular}
%\centering
%\caption{Iteration counts for highly advection-dominated case.}
%\begin{tabular}{| c | c | c | c | c |}
%\hline
%Methods/Models & $20^2$ & $40^2$ & $80^2$ & $160^2$ \\
%\hline
%%Time(s) & 149.8 & 225.3 & 238 \\
%AMG & - & - & - & -\\
%CPR-AMG(1) & 15.0 & 15.3 & 15.6 & 16.6 \\
%CPR-AMG(2) & 28 & 29.7 & 30.6 & 32.0 \\
%BF & 10.7 & 10.8 & 11.3 & 12.3 \\
%\hline
%\end{tabular}
%\end{table}
In this example, we compare the performance of the different strategies for a problem in which gravity plays a dominant role. The domain is a square box of size $20 \times 20$ meters. The absolute permeability is a homogeneous field of $100$ millidarcy. 
Water is injected into the domain through the boundary at the top left corner, and the outlet is at the top right corner. The rate of injection is 5 $m^3$/day. 
For spatial discretization, we use uniform grids of size $20 \times 20$, $40 \times 40$, $80 \times 80$, and $160 \times 160$ respectively.
The initial conditions are the same as the heterogeneous two-dimensional example above. The time steps are 10, 4, and 1 days, and the final times are 20, 8, and 2 days respectively.
\begin{table}
\centering
\caption{Iteration counts for diffusion-dominated case with gravity, time step $dt = 10$ days.} \label{2d_grav_diff}
\begin{tabular}{| c | c | c | c | c |}
\hline
Methods/Mesh sizes & $20^2$ & $40^2$ & $80^2$ & $160^2$\\
\hline
%Time(s) & 149.8 & 225.3 & 238 \\
AMG & 7 & 7 & 7 & 7\\
CPR-AMG(1) & 15.1 & 25.9 & 49.4 & 95.2 \\
CPR-AMG(2) & 22.0 & 30.9 & 38.1 & 40.7 \\
BF & 19.9 & 21.0 & 21.1 & 21.1 \\
\hline
\end{tabular}
\centering
\caption{Iteration counts for advection-dominated case with gravity, time step $dt = 4$ days.} \label{2d_grav_adv}
\begin{tabular}{| c | c | c | c | c |}
\hline
Methods/Mesh sizes & $20^2$ & $40^2$ & $80^2$ & $160^2$\\
\hline
%Time(s) & 149.8 & 225.3 & 238 \\
AMG & - & - & - & 17.3 \\
CPR-AMG(1) & 17.9 & 17.8 & 16.1 & 25.2 \\
CPR-AMG(2) & 30.3 & 29.8 & 22.7 & 30.2 \\
BF & 13.0 & 17.0 & 21.1 &  23.9 \\
\hline
\end{tabular}
\centering
\caption{Iteration counts for highly advection-dominated case with gravity, time step $dt = 1$ day.} \label{2d_grav_hyp}
\begin{tabular}{| c | c | c | c | c |}
\hline
Methods/Mesh sizes & $20^2$ & $40^2$ & $80^2$ & $160^2$ \\
\hline
%Time(s) & 149.8 & 225.3 & 238 \\
AMG & - & - & - & -\\
CPR-AMG(1) & 18.6 & 19.6 & 19.7 & 18.8 \\
CPR-AMG(2) & 31.5 & 34.6 & 36.6 & 36.3 \\
BF & 13.1 & 9.3 & 12.0 & 16.4 \\
\hline
\end{tabular}
\end{table}

The diffusion-dominated case, shown in Table \ref{2d_grav_diff}, 
%we observe 
exhibits the same pattern as in the previous example: the AMG preconditioner is the most efficient method, followed by the block factorization method, CPR-AMG(2), and CPR-AMG(1). 
AMG and the block factorization method exhibit optimal 
%scaling 
performance
with respect to problem size. The number of iterations for CPR-AMG(2) also seems to reach a plateau 
%at the largest mesh size. 
as the mesh size is refined.
%Meanwhile, 
In contrast, the performance of CPR-AMG(1) does not %show good scaling result 
scale well with respect to mesh size for this case, taking about twice the number of iterations for each level of mesh refinement.

The results for the advection-dominated case are shown in Table \ref{2d_grav_adv}. 
%lacks robustness 
The AMG method is not robust
and only converges for the largest mesh size (for which it takes the fewest iterations).
The block factorization preconditioner %exhibits strong robustness and scaling properties, 
is highly robust and also appears to require iteration counts tending to a constant as
the mesh is refined.
%as the number of iterations tends to get flat at the largest mesh size.  
The performance of CPR-AMG(2) is consistent except for the $80 \times 80$ mesh. Although it requires more iterations than CPR-AMG(1), this method shows promising scaling property, similar to the previous example, since the number of iterations does not grow as the mesh is refined. CPR-AMG(1) performs quite well for this case, but it still exhibits poor scalability as the number of iterations grow quickly between $80 \times 80$ and $160 \times 160$.

In the strongly advection-dominated case, AMG diverges for all mesh sizes. The new block factorization is the most efficient method in this case, requiring the smallest number of iterations across all mesh sizes. Here, CPR-AMG(1) is more efficient than CPR-AMG(2), requiring about half the number of iterations. Both CPR-AMG(1) and CPR-AMG(2) show good scaling property in this case. The scaling result for the block factorization method is not as clear as in the diffusion and advection dominated cases, but 
%it is possible that the mesh is not large 
we suspect that the mesh is not fine enough for a consistent pattern to emerge.

Besides varying the mesh size, we also experimented with changing the time step size for a fixed mesh of $80 \times 80$ for the same problem. The results are reported in Table \ref{varying_time_steps}. The final time for the simulation is 8 days. The results are reported in Table \ref{varying_time_steps}.
%\begin{table}[H]
%\centering
%\begin{tabular}{| c | c | c | c | c | c | c | c | c | c | c | c | c |}
%\hline
%\multirow{2}{*}{Methods/Models} & \multicolumn{3}{c|}{1.25 days} & \multicolumn{3}{c|}{2.5 days} & \multicolumn{3}{c|}{5 days} & \multicolumn{3}{c|}{10 days}\\
%\cline{2-13}
%& NI & LI & LI/NI & NI & LI & LI/NI & NI & LI & LI/NI & NI & LI & LI/NI \\
%\hline
%CPR-AMG(1) & 54 & 817 & 15.1 & 65 & 987 & 15.2 & 76 & 1213 & 16.0 & 75 & 1562 & 20.8 \\
%CPR-AMG(2) & 54 & 1385 & 25.6 & 65 & 1396 & 21.5 & 76 & 1846 & 24.3 & 75 & 2983 & 39.8 \\
%BF & 54 & 812 & 15.0 & 65 & 1151 & 17.7 & 76 & 1562 & 20.6 & 75 & 1900 & 25.3 \\
%\hline
%\end{tabular}
%\caption{Results for the advection-dominated case $P_0 = 10^4$.}
%\end{table}

%\begin{table}
%\centering
%\begin{tabular}{| c | c | c | c | c | c | c | c | c | c | c | c | c |}
%\hline
%\multirow{2}{*}{Methods/Time steps} & \multicolumn{3}{c|}{$dt$ = 1 day} & \multicolumn{3}{c|}{$dt$ = 2 days} & \multicolumn{3}{c|}{$dt$ = 4 days} & \multicolumn{3}{c|}{$dt$ = 8 days}\\
%\cline{2-13}
%& NI/TS & LI & LI/NI & NI/TS & LI & LI/NI & NI/TS & LI & LI/NI & NI & LI & LI/NI \\
%\hline
%CPR-AMG(1) & 99 & 1672 & 16.9 & 68 & 1083 & 15.9 & 46 & 738 & 16.1 & 28 & 544 & 19.4 \\
%CPR-AMG(2) & 99 & 2870 & 29.0 & 68 & 1630 & 24.0 & 46 & 1042 & 22.7 & 28 & 761 & 27.2 \\
%BF & 99 & 1657 & 16.7 & 68 & 1292 & 19.0 & 46 & 968 & 21.1 & 28 & 645 & 23.0 \\
%\hline
%\end{tabular}
%\caption{Results for the advection-dominated case with gravity $P_0 = 10^4$.} \label{varying_time_steps}
%\end{table}
\begin{table}
\captionsetup{justification=centering}
\caption{Results for the advection-dominated case with gravity $P_0 = 10^4$. \\ NI/TS is the number of Newton iteration per time step.} \label{varying_time_steps}
\centering
\begin{tabular}{| c | c | c | c | c | c | c | c | c |}
\hline
\multirow{2}{*}{Methods/Time steps} & \multicolumn{2}{c|}{$dt$ = 1 day} & \multicolumn{2}{c|}{$dt$ = 2 days} & \multicolumn{2}{c|}{$dt$ = 4 days} & \multicolumn{2}{c|}{$dt$ = 8 days}\\
\cline{2-9}
& NI/TS & LI/NI & NI/TS & LI/NI & NI/TS & LI/NI & NI/TS & LI/NI \\
\hline
CPR-AMG(1) & 12.4 & 16.9 & 17 & 15.9 & 23 & 16.1 & 28 & 19.4 \\
CPR-AMG(2) & 12.4 & 29.0 & 17 & 24.0 & 23 & 22.7 & 28 & 27.2 \\
BF & 12.4 &16.7 & 17 & 19.0 & 23 & 21.1 & 28 & 23.0 \\
\hline
\end{tabular}
\end{table}
Since AMG does not converge in this experiment, we exclude it from the results. From Table \ref{varying_time_steps}, it is clear that as the time step gets larger, Newton's method takes more iterations to converge. For $dt=$ 8 days, there is only one time step and it is the most difficult case. CPR-AMG(1) number of iterations is not significantly affected by the time step except for the largest time step size of 8 days. Meanwhile, CPR-AMG(2) number of iterations decreases as the time step gets larger, but goes up again at $dt = 8$ days. The block factorization method shows consistent increase in the number of iterations for larger time steps. Overall, there is not much of a difference in terms of iteration counts for these three methods, but it is worth noting that the block factorization method is much faster than the others in terms of run time, as it does not require a global ILU solve.

\subsection{Behavior of Eigenvalues}
It is often possible to obtain insight into the properties of
preconditioning operators from the eigenvalues of the preconditioned
matrix $JM^{-1}$.
In particular, recall a standard analysis of the convergence behavior of 
GMRES for solving the preconditioned system (\ref{precon-system}) 
\cite{Saad03}.
Assume the preconditioned matrix is diagonalizable, 
$JM^{-1} = V \Lambda V^{-1}$,
where $\Lambda$ is a diagonal matrix containing the eigenvalues of the
preconditioned matrix and the columns of $V$ are the corresponding eigenvectors.
If $c_k=M^{-1}{\hat c}_k$ are the iterates obtained at the $k$th step of
GMRES iteration, with  residual $r_k = q - J c_k$, then
\begin{equation} \label{GMRES-bound}
\frac{\|r_k\|_2}{\|r_0\|_2} \le \|V\|_2 \, \|V^{-1}\|_2
\min_{p_k(0)=1} \max_{\lambda \in \sigma(JM^{-1})} |p_k(\lambda)|,
\end{equation}
where the minimum in (\ref{GMRES-bound}) is over all polynomials of degree
at most $k$ that have the value $1$ at the origin,
$\sigma(JM^{-1})$ is the set of eigenvalues of $JM^{-1}$,
and the norm is the vector Euclidian norm.
Thus, a good preconditioner tends to produce a preconditioned operator with
a compressed spectrum whose entries are not near the origin.
In this section, we explore the behavior of the eigenvalues of the
preconditioned matrix with an eye toward understanding the effects of 
features of the discrete problem such as
discretization mesh size and
qualitative features of the model such as the relative weights of 
diffusion and advection and the degree of coupling between the components.

\begin{figure}
\begin{picture}(400,480)
\put(10,318) {\includegraphics[width=8.0cm]{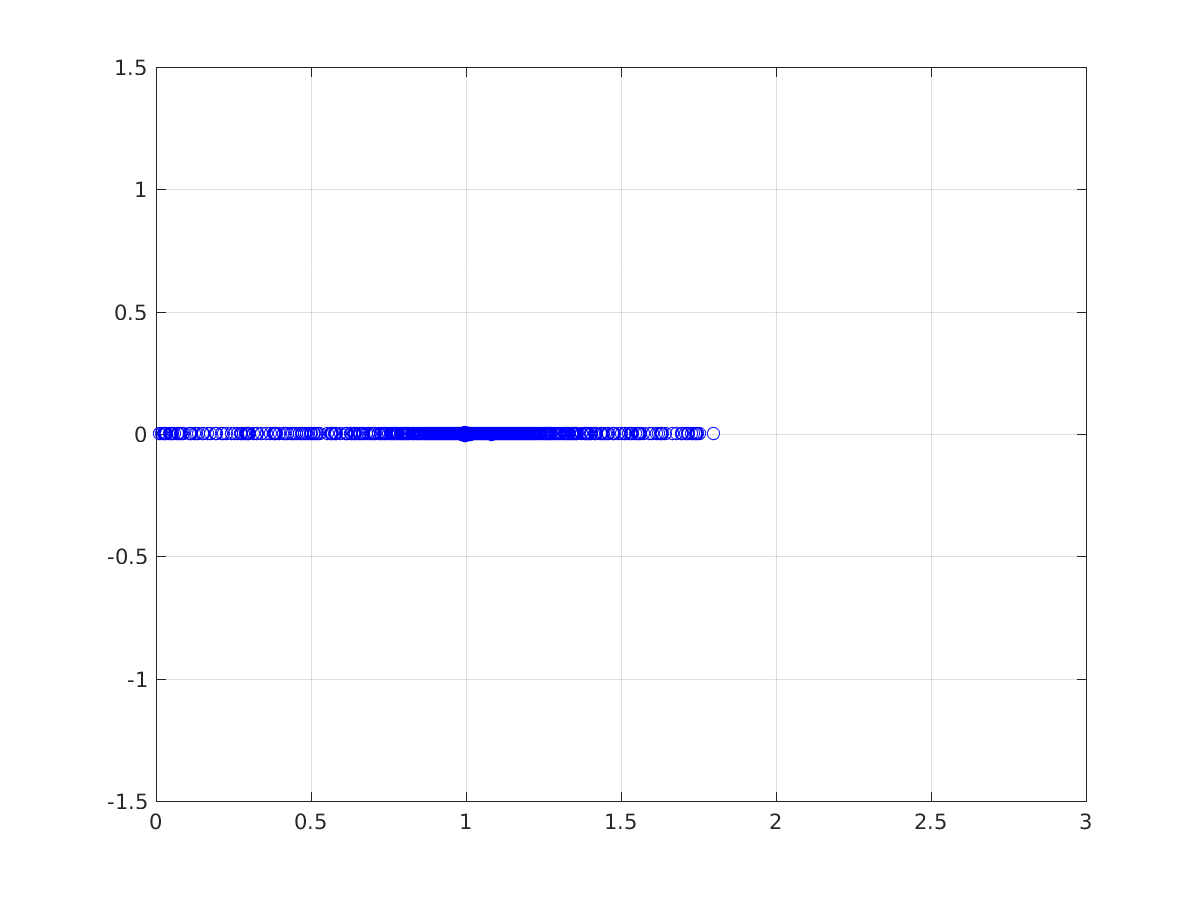}}
\put(125,455){Diffusion-dominated} \put(125,440){CPR-AMG(1)}
\put(220,318){\includegraphics[width=8.0cm]{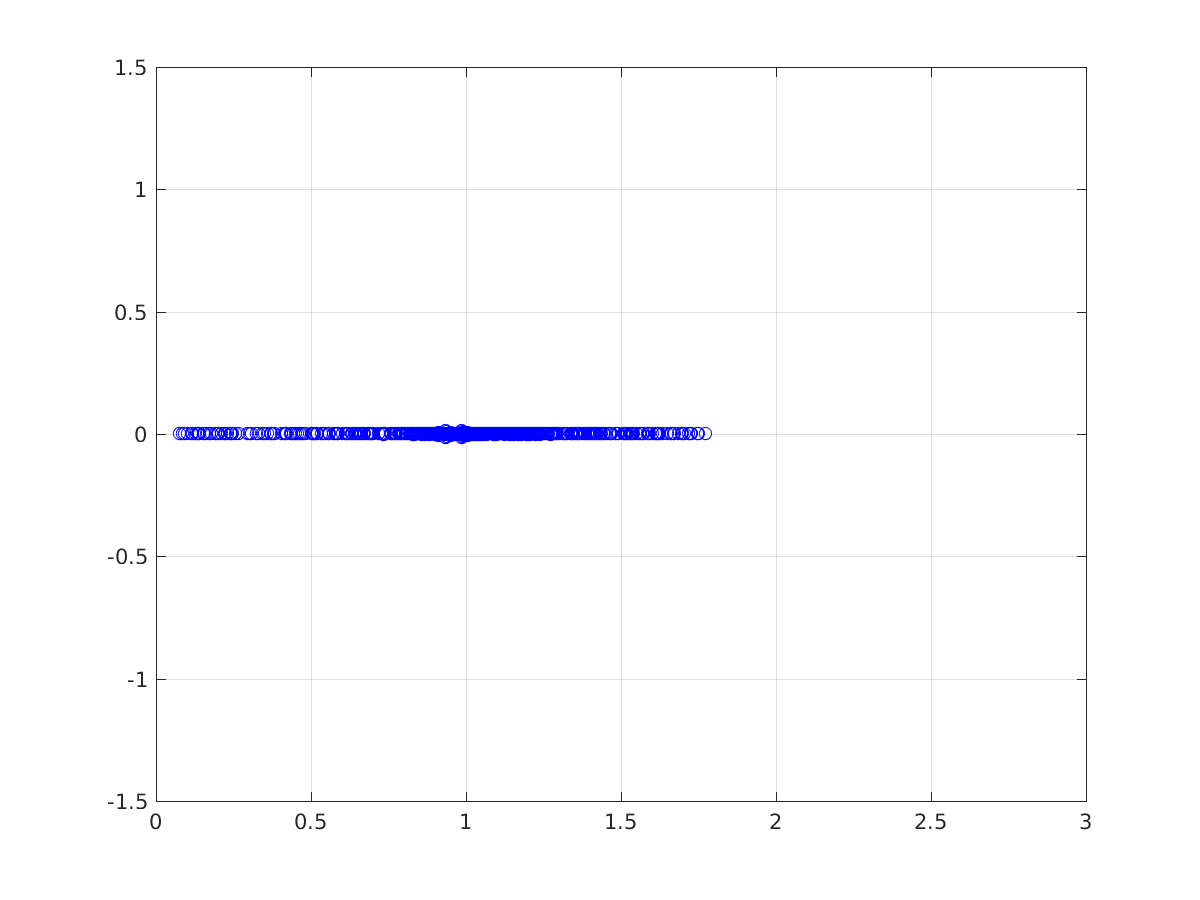}}
\put(328,455){Advection-dominated} \put(328,440){CPR-AMG(1)}
\put(10,159) {\includegraphics[width=8.0cm]{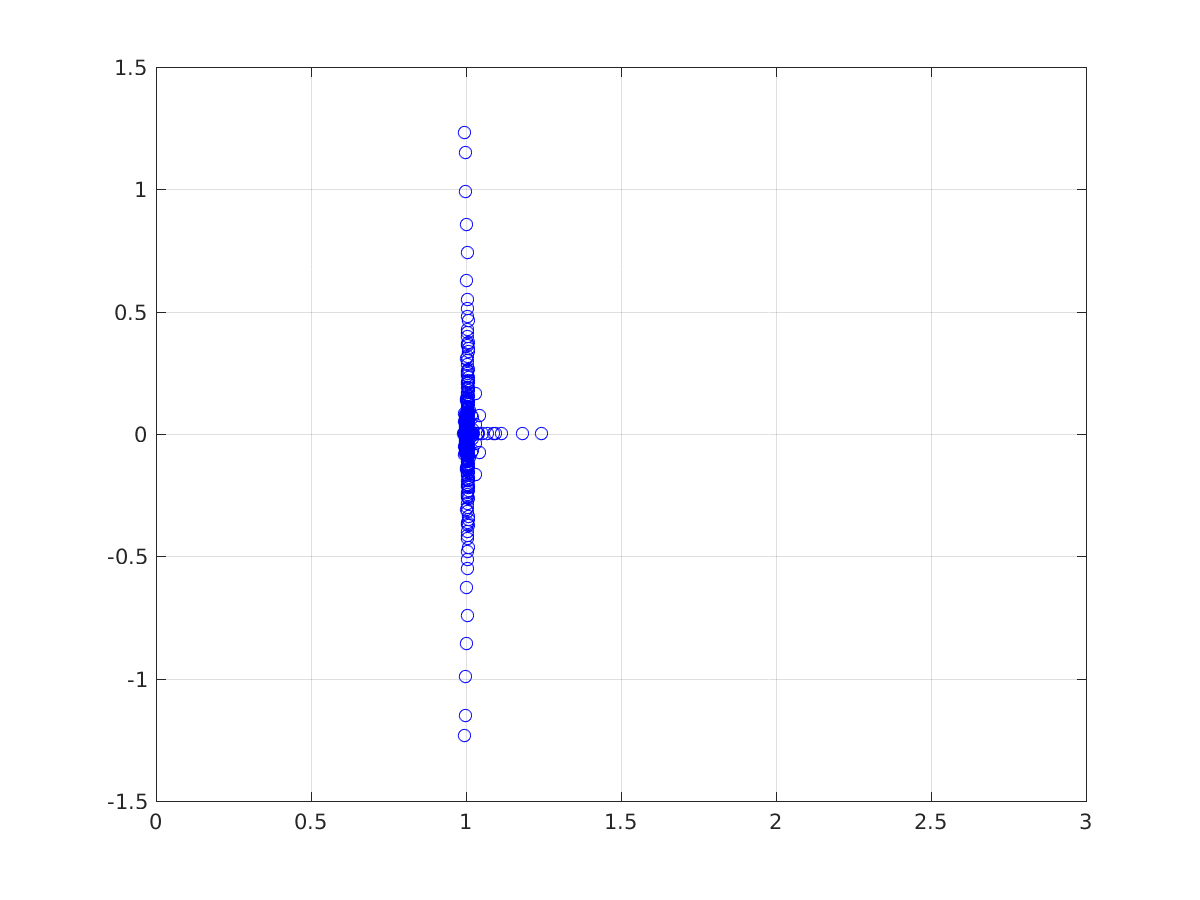}}
\put(125,298){Diffusion-dominated} \put(125,283){CPR-AMG(2)}
\put(220,159){\includegraphics[width=8.0cm]{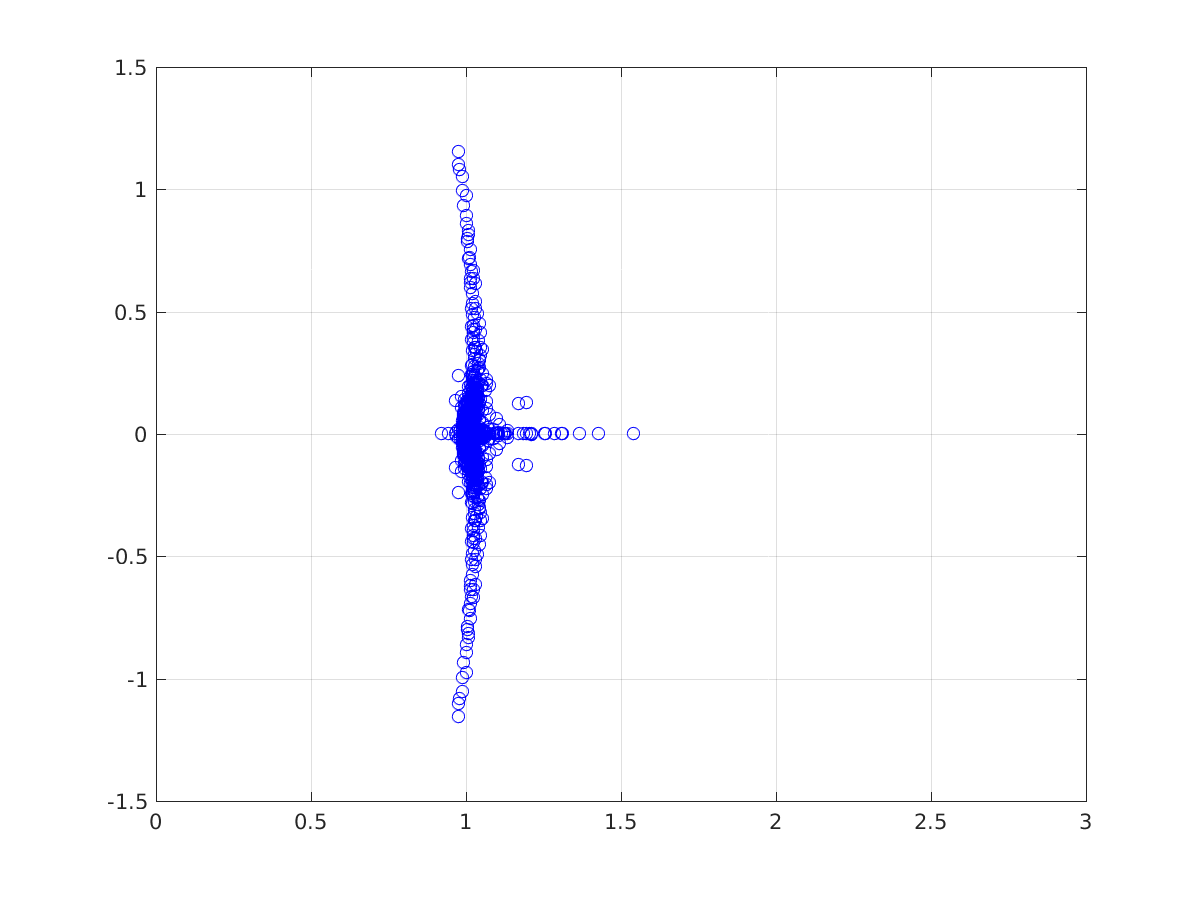}}
\put(328,298){Advection-dominated} \put(328,283){CPR-AMG(2)}
\put(10,0)   {\includegraphics[width=8.0cm]{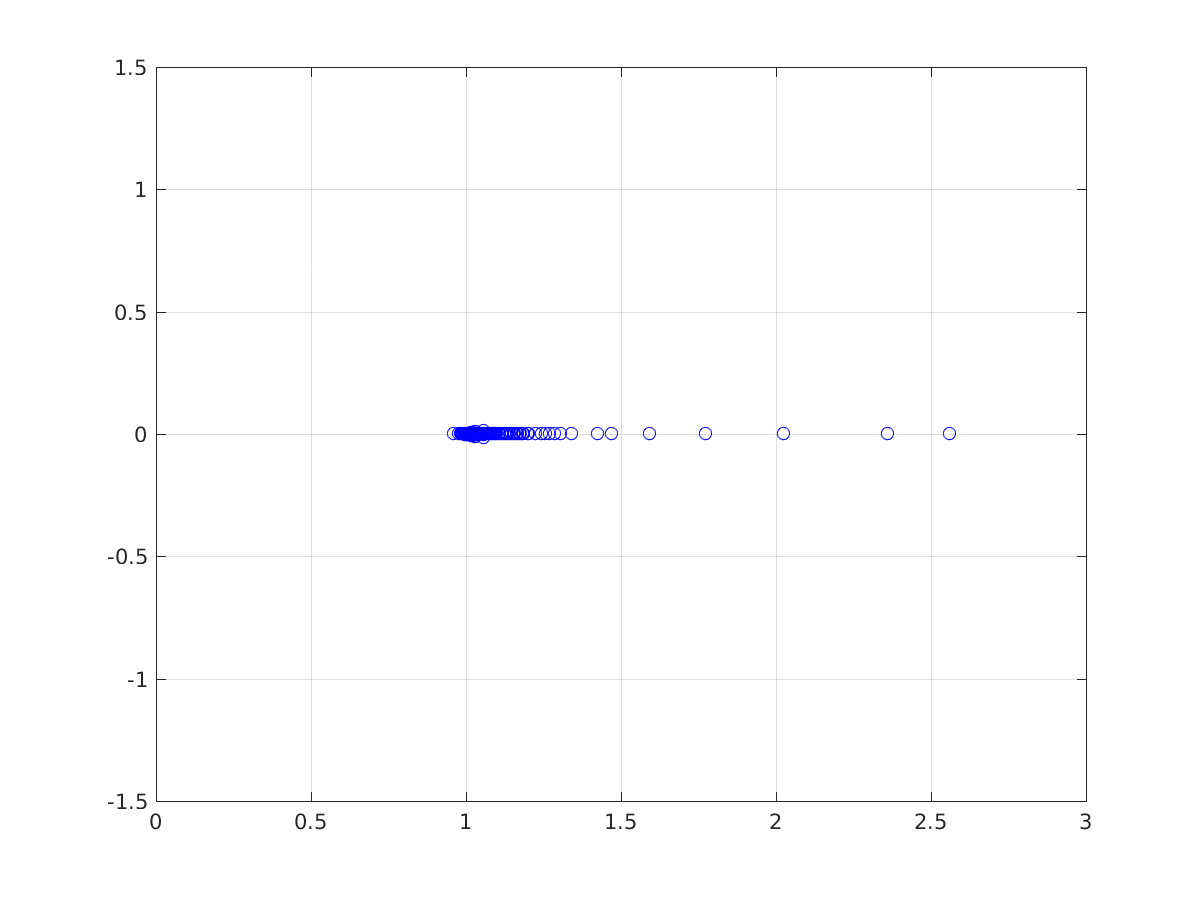}}
\put(125,141){Diffusion-dominated} \put(125,126){BF}
\put(220,0)  {\includegraphics[width=8.0cm]{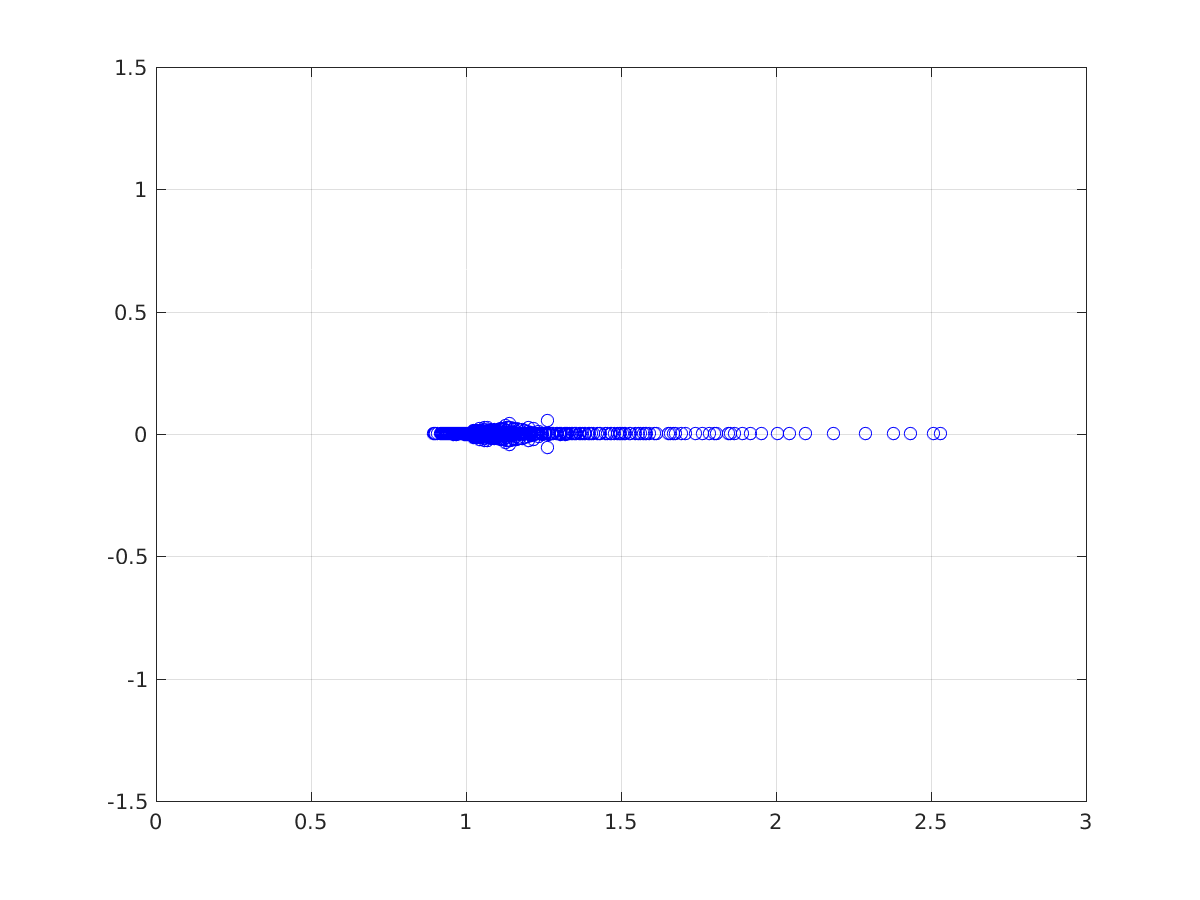}}
\put(328,141){Advection-dominated} \put(328,126){BF}
\end{picture}
%\centering
%\subfloat[ILU]{\includegraphics[width=0.45\textwidth]{evals_diff_1e5_ilu_40x40.png}}
%\subfloat[CRP-AMG(1)]{\includegraphics[width=0.45\textwidth]{evals_diff_1e5_cpr1_40x40.png}}
%\subfloat[CPR-AMG(2)]{\includegraphics[width=0.45\textwidth]{evals_diff_1e5_cpr2_40x40.png}}
%\subfloat[Block Factorization]{\includegraphics[width=0.45\textwidth]{evals_diff_1e5_bf_40x40.png}} \\
\caption{Eigenvalues of preconditioned systems for different strategies, applied to
  the diffusion-dominated Example 1 (left) and advection-dominated Example 2 (right).}
\label{evalues_plot}
\end{figure}
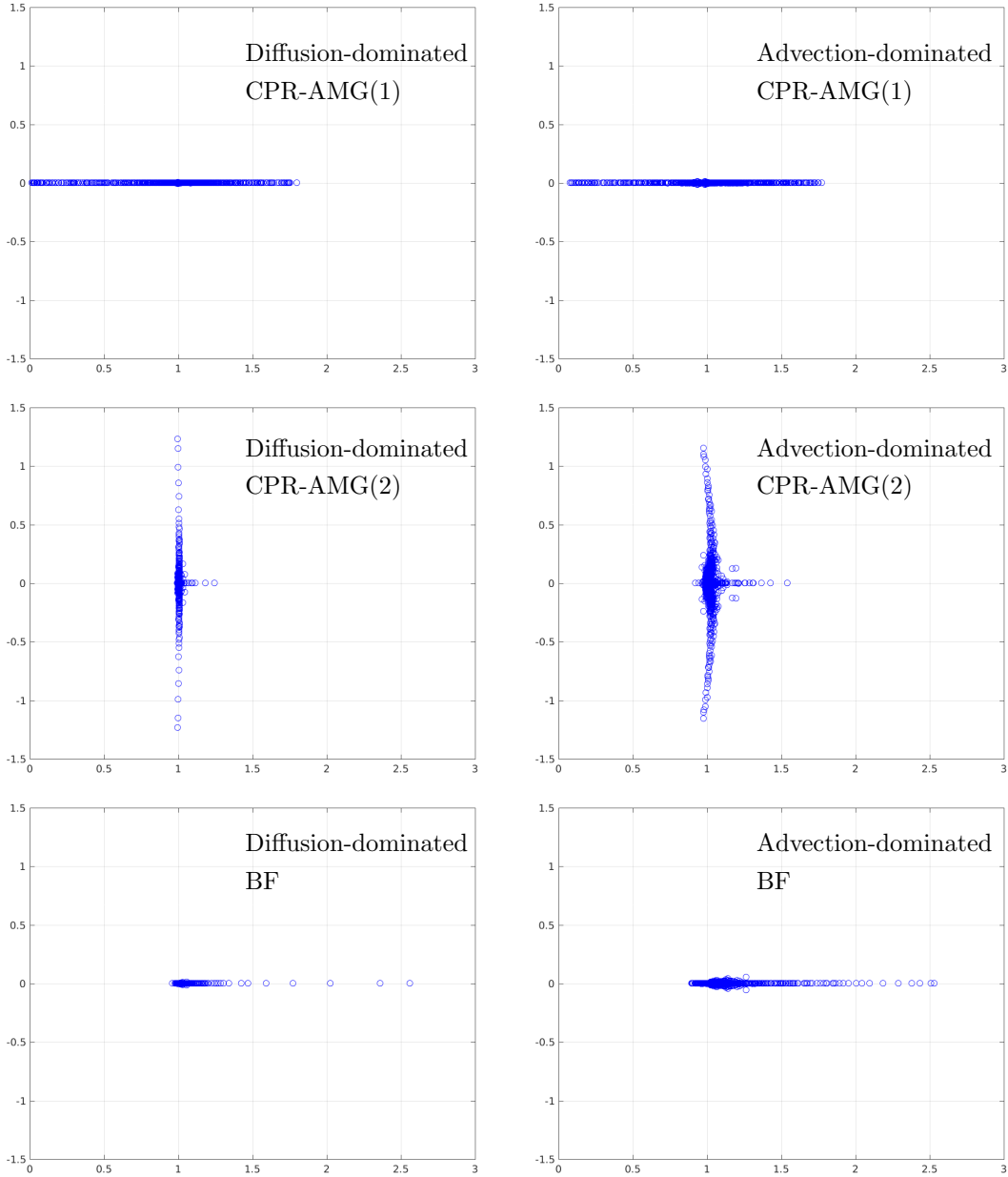

Figure \ref{evalues_plot} gives a representative depiction of the eigenvalues
of preconditioned operators for three of the preconditioners considered.
These results are for benchmark problems for which performance is considered
in section \ref{2D-oil-water}, the two-dimensional linear oil-water model discretized
on a $100 \times 20$ grid.
The plots on the left side of the figure show eigenvalues for the
diffusion-dominated case (for which solution performance is shown in
Table \ref{performance_q5_1}),
and those on the right show eigenvalues for the advection-dominated case
(performance in Table \ref{performance_q5_2}). \footnote{These computations were done using the {\tt eig} function in Matlab, and they use Matlab backslash to perform the actions of the inverses of $A_{11}$, $A_{22}$ and the modified Schur complement. This contrasts with the solution algorithms tested, which approximate these operations using one AMG V-cycle.}

These displays indicate that the spectra for the preconditioned systems
for the CPR-AMG(2) and BF preconditioners are bounded away from the origin,
whereas for the CPR-AMG(1) preconditioner there are many small eigenvalues.
Performance of CPR-AMG(1) improves in the advection-dominated case, and
the smallest associated eigenvalues are somewhat further from the origin.
In contrast, the latter two preconditioners are largely unchanged in
the advection-dominated case, where they are still effective, and the associated eigenvalues are also contained in similarly structured regions far from the origin.
We believe the superior performance of the BF preconditioner comes from its
greater emphasis on the coupling between pressure and saturation, derived
from use of the approximate Schur complement (\ref{SIMPLE-Schur}).

%In order to understand the effectiveness of the two CPR-AMG approaches, we compute all the eigenvalues of the preconditioned operator by solving the generalized eigenvalue problem
%\begin{align}
%Av = \lambda Mv,
%\end{align}
%where $M$ is the preconditioner whose inverse is defined in \eqref{m_comb}, \eqref{m_add}, and \eqref{m_bf}. We also compare the results to the eigenvalues of $M$ constructed using ILU, which we know is not a good preconditioner for these types of problems. We use a mesh of size $40 \times 40$. The eigenvalues of four strategies

%Things that can be included
%\begin{itemize}
%\item The same figures for the advection-dominated case.
%\item The same figures with overlapping Ritz values.
%\end{itemize}

\subsection{Three-dimensional Problem}
We use a homogeneous permeability field of $100$ millidarcy, and the grid is stretched to induce anisotropy. The model dimensions are $25 \times 100 \times 6$ meters and the cell size is $0.5\times 1 \times 0.05$ meter. Thus, the mesh is $50\times 100 \times 120$, and the problem has 1.2 million unknowns in total. Water is injected into the domain at one bottom corner and the outlet is at the opposite corner. The injection rate is 0.75 $m^3/day$. The parameters for the capillary pressure model is from example 1 of Table \ref{pc_params}. The simulation is run for 100 days with time step $\Delta t = 20$ days.
\begin{table}
\centering
\captionsetup{justification=centering}
\caption{Performance in the 3D case for the set of parameters in example 1 of Table \ref{pc_params}} \label{performance_3d_1}
\begin{tabular}{| c | c | c | c | c | c | c | c | c |}
\hline
\multirow{2}{*}{Methods/Models} & \multicolumn{4}{c|}{Linear} & \multicolumn{4}{c|}{Brooks Corey} \\
\cline{2-9}
& NI & LI & LI/NI & Time & NI & LI & LI/NI & Time \\
\hline
AMG & 16 & 282 & 17.6 & 103.1 & 20 & 452 & 22.6 & 144.7\\
CPR-AMG(1) & 16 & 2698 & 168.6 & 803.2 & 20 & 6069 & 303.45 & 1940.8\\
CPR-AMG(2) & 16 & 712 & 44.5 & 299.5 & 20 & 1900 & 95.0 & 741.1 \\
BF & 16 & 355 & 22.2 & 133.6 & 20 & 752 & 37.6 & 231.1\\
\hline
\end{tabular}
\end{table}
Table \ref{performance_3d_1} shows the performance results of the diffusion-dominated case for this 3D example, which are consistent with those of the previous two-dimensional example. AMG preconditioner shows the best results for both the iteration counts per Newton step and the time it takes to complete the simulation for both capillary pressure models. CPR-AMG(2) does not perform quite as well as AMG, but it is much more efficient than CPR-AMG(1) for both performance measures and capillary pressure models. As in the two-dimensional case, the new block factorization method performs well, requiring fewer than half the iterations than CPR-AMG(2) for both the linear and Brooks-Corey models, and running in about one third the CPU time.

We also tested the three-dimensional SPE10 problem with the linear model of capillary pressure for the different preconditioning strategies. Here, AMG diverges even for the diffusion-dominated case, even though it was the most efficient method for the two-dimensional example. The block factorization method is about four times faster than CPR-AMG(2) and five times faster than CPR-AMG(1) in the diffusion dominated case. CPR-AMG(2) still outperforms CPR-AMG(1) both in terms of iteration counts and run time, but the margin is smaller than for the two-dimensional problem. In the advection-dominated case, unlike in the two-dimensional example, CPR-AMG(1) is more efficient than CPR-AMG(2), requiring about 23\% fewer number iterations and 45\% run time. The block factorization approach is still the most efficient method, taking fewer than half the number of iterations and less than half the run time of CPR-AMG(1). We also note that the number of iterations for the block factorization method is very consistent with respect to the characteristics of the problem, i.e. it does not change significantly whether the problem is diffusion-dominated, advection-dominated, or strongly advection-dominated.
\begin{table}
\centering
\begin{tabular}{| c | c | c | c | c |}
\hline
\multirow{2}{*}{Methods/Models} & \multicolumn{4}{c|}{Linear} \\
\cline{2-5}
& NI & LI & LI/NI & Time (s)\\
\hline
AMG & - & - & - & - \\
CPR-AMG(1) & 17 & 2410 & 141.8 & 614.14 \\
CPR-AMG(2) & 17 & 1661 & 97.7 & 448.21 \\
BF & 17 & 490 & 28.8 & 121.71 \\
\hline
\end{tabular}
\captionsetup{justification=centering}
\caption{Performance of the four methods for linear entry pressure $P_0 = 10^6$.}
\begin{tabular}{| c | c | c | c | c |}
\hline
\multirow{2}{*}{Methods/Models} & \multicolumn{4}{c|}{Linear} \\
\cline{2-5}
& NI & LI & LI/NI & Time (s)\\
\hline
AMG & - & - & - & - \\
CPR-AMG(1) & 18 & 1122 & 62.3 & 354.38 \\
CPR-AMG(2) & 18 & 1554 & 86.3 & 657.12 \\
BF & 18 & 474 & 26.3 & 157.24 \\
\hline
\end{tabular}
\captionsetup{justification=centering}
\caption{Performance of the four methods for linear entry pressure $P_0 = 10^5$.}
\end{table}

\subsection{Scaling Results}
To perform a scalability study, we run a test problem on a box of dimensions $20\times 20\times 20$ meters. The initial mesh is $20\times 20\times 20$ and is repeatedly refined in the z-direction. The domain has constant material properties. The parameters for the water retention models are listed in example 4 of Table \ref{pc_params}. Note that this set of parameters corresponds to a diffusion-dominated problem.
\begin{figure}
\centering
\includegraphics[width=0.45\textwidth]{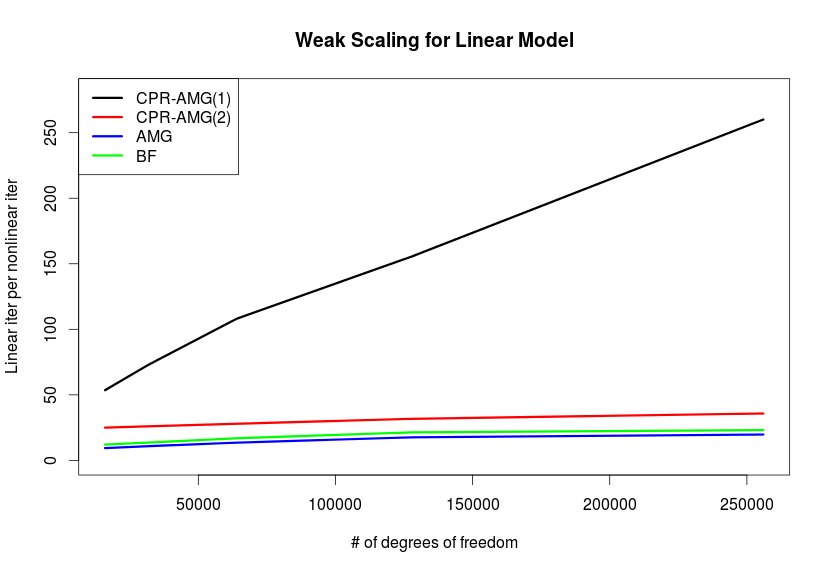}
\includegraphics[width=0.45\textwidth]{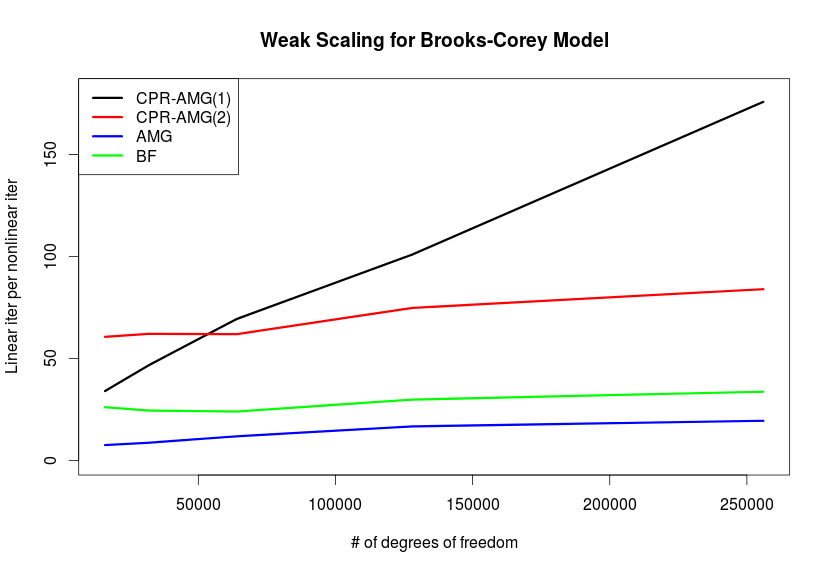}
\caption{Weak scaling for different strategies.} \label{weak_scaling}
\end{figure}
The results shown in Figure \ref{weak_scaling}. indicate that the performance of the block factorization, CPR-AMG(2), and AMG methods is independent of the mesh size. The number of linear iterations per Newton step does not grow as the mesh is refined which is optimal multigrid performance. The block factorization method's performance is nearly identical to that of AMG for the linear model, and still quite close for the Brooks-Corey model, compared to CPR-AMG(2). CPR-AMG(1), however, does not scale as well as the other two methods. The linear iteration counts for CPR-AMG(1) grows linearly as the mesh is refined.
%\begin{figure}[H]
%\centering
%\includegraphics[width=0.45\textwidth]{strong_bc.png} 
%\includegraphics[width=0.45\textwidth]{strong_iter.png} 
%\caption{Strong scaling for different strategies.} \label{timing_bc}
%\end{figure}
%For a strong scaling study, the mesh is fixed at $40\times 40\times 80$ with 256000 degrees of freedom. We run the problem on 16, 32, 64, 128 and 256 cores. Timings is shown in Figure \ref{timing_bc}.

\section{Conclusions}
In this work, we have implemented a fully implicit parallel isothermal two-phase flow simulator along with four different preconditioning strategies to solve the linear systems resulting from linearization of the coupled equations, and we have tested the performance of these methods as preconditioners for GMRES. We have also developed a new block factorization preconditioner whose performance is robust and efficient across all benchmark problems studied. In contrast, although AMG preconditioning applied to the coupled systems is the most efficient choice in some cases (both two-dimensional and three-dimensional diffusion-dominated examples), it exhibits slow convergence and sometimes diverges for advection-dominated cases. The new block factorization preconditioner achieves consistently low iteration counts across all the tests and varying examples of capillary pressure, and it scales optimally with problem size. The combinative CPR-AMG(1), though robust across all the tests, is the least efficient method, with the exception of the near hyperbolic case where it is faster than CPR-AMG(2). The additive CPR-AMG(2) method performs well in most cases except the strongly advection-dominated case. It also scales optimally with problem size for both advection-dominated and diffusion-dominated case. %The new block factorization is very efficient and robust, achieving consistently low number of iteration counts across all the tests and different cases for capillary pressure models. Except for the diffusion-dominated case, the new block factorization method outperforms other methods in terms of both iteration counts and run time. It also scales optimally with problem size.

\section*{Acknowledgement}
This work was supported by the US Department of Energy under grant DE-SC0009301, the U.S. National Science Foundation under grant DMS1418754, and by the Department of Energy at Los Alamos National Laboratory under contract DE-AC52-06NA25396 through the DOE Advanced Simulation Capability for Environmental Management (ASCEM) program, and the LANL Summer Student Program at the Center for Nonlinear Studies (CNLS).

\newpage
\bibliographystyle{plain}
\bibliography{siam_si.bib}

\begin{thebibliography}{10}

\bibitem{Amanzi}
ASCEM.
\newblock \url{http://esd.lbl.gov/research/projects/ascem/thrusts/hpc/}, 2009.

\bibitem{Axelsson94}
O.~Axelsson.
\newblock {\em Iterative Solution Methods}.
\newblock Cambridge University Press, 1994.

\bibitem{Aziz79}
K.~Aziz and A.~Settari.
\newblock {\em Petroleum Reservoir Simulation}.
\newblock Applied Science Publishers, 1979.

\bibitem{Bastian99}
P.~Bastian.
\newblock Numerical computation of multiphase flow in porous media.
\newblock Habilitationsschrift, 1999.

\bibitem{BastianHelmig99}
P.~Bastian and R.~Helmig.
\newblock Efficient fully-coupled solution techniques for two-phase flow in
  porous media: Parallel multigrid solution and large scale computations.
\newblock {\em Advances in Water Resources}, 23(3):199 -- 216, 1999.

\bibitem{Behie82}
A.~Behie and P.K.W. Vinsome.
\newblock Block iterative methods for fully implicit reservoir simulation.
\newblock {\em Society of Petroleum Engineers Journal}, 22(05):658--668, oct
  1982.

\bibitem{BrooksCorey64}
R.~Brooks and A.~Corey.
\newblock {\em {Hydraulic Properties of Porous Media}}.
\newblock Hydrology Papers, Colorado State Univeristy, 1964.

\bibitem{Christie01}
M.A. Christie and M.J. Blunt.
\newblock Tenth {SPE} comparative solution project: A comparison of upscaling
  techniques.
\newblock In {\em {SPE} Reservoir Simulation Symposium}. Society of Petroleum
  Engineers ({SPE}), 2001.

\bibitem{Cleary98}
A.J. Cleary, R.D. Falgout, V.E. Henson, and J.E. Jones.
\newblock Coarse-grid selection for parallel algebraic multigrid.
\newblock In Alfonso Ferreira, José Rolim, Horst Simon, and Shang-Hua Teng,
  editors, {\em Solving Irregularly Structured Problems in Parallel}, volume
  1457 of {\em Lecture Notes in Computer Science}, pages 104--115. Springer
  Berlin Heidelberg, 1998.

\bibitem{Clees10}
T.~Clees and L.~Ganzer.
\newblock An efficient algebraic multigrid solver strategy for adaptive
  implicit methods in oil-reservoir simulation.
\newblock {\em {SPE} Journal}, 15(03):670--681, 2010.

\bibitem{ECoon_JDMoulton_SPainter_2014a}
E.~Coon, J.~D. Moulton, and S.~Painter.
\newblock Managing complexity in simulations of land surface and near-surface
  processes.
\newblock Technical Report {LA-UR} 14-25386, Applied Mathematics and Plasma
  Physics Group, Los Alamos National Laboratory, 2014.

\bibitem{Corey94}
A.~T. Corey.
\newblock {\em Mechanics of Immiscible Fluids in Porous Media}.
\newblock Water Resources Pubns, 1995.

\bibitem{Dawson97}
C.N. Dawson, H.~Klie, M.F. Wheeler, and C.S. Woodward.
\newblock A parallel, implicit, {cell-centered} method for {two-phase} flow
  with a preconditioned {Newton–Krylov} solver.
\newblock {\em Computational Geosciences}, 1(3-4):215--249, 1997.

\bibitem{Falgout06}
R.D. Falgout, J.E. Jones, and U.M. Yang.
\newblock The design and implementation of hypre, a library of parallel high
  performance preconditioners.
\newblock In AreMagnus Bruaset and Aslak Tveito, editors, {\em Numerical
  Solution of Partial Differential Equations on Parallel Computers}, volume~51
  of {\em Lecture Notes in Computational Science and Engineering}, pages
  267--294. Springer Berlin Heidelberg, 2006.

\bibitem{Falgout02}
R.D. Falgout and U.M. Yang.
\newblock {HYPRE}: a library of high performance preconditioners.
\newblock In {\em Preconditioners,” Lecture Notes in Computer Science}, pages
  632--641, 2002.

\bibitem{Henson00}
V.E. Henson and U.M. Yang.
\newblock {BoomerAMG}: a parallel algebraic multigrid solver and
  preconditioner.
\newblock {\em Applied Numerical Mathematics}, 41:155--177, 2000.

\bibitem{Klie96}
H.~Klie, M.~Rame, and M.F. Wheeler.
\newblock Two-stage preconditioners for inexact {Newton} methods in multiphase
  reservoir simulation.
\newblock Technical report, Center for Research on Parallel Computation, Rice
  University, 1996.

\bibitem{Lacroix03}
S.~Lacroix, Yu. Vassilevski, J.~Wheeler, and M.F. Wheeler.
\newblock Iterative solution methods for modeling multiphase flow in porous
  media fully implicitly.
\newblock {\em {SIAM} J. Sci. Comput.}, 25(3):905--926, jan 2003.

\bibitem{Lu08}
B.~Lu.
\newblock {\em Iteratively Coupled Reservoir Simulation for Multiphase Flow in
  Porous Media}.
\newblock PhD thesis, Center for Subsurface Modeling, {ICES}, The University of
  Texas - Austin, 2008.

\bibitem{Masson12}
R.~Masson, R.~Scheichl, P.~Quandalle, C.~Guichard, R.~Eymard, R.~Herbin, and
  K.~Brenner et~al.
\newblock {CPR-AMG} preconditioner for multiphase porous media flows: a few
  experiments.
\newblock {La page num\'erique du LATP - UMR 7353- Marseille}, 2012.

\bibitem{Osei-Kuffuor15}
D.~Osei-Kuffuor, L.~Wang, R.~D. Falgout, and I.~D. Mishev.
\newblock An algebraic multigrid solver for fully-implicit solution methods in
  reservoir simulation.
\newblock SIAM Geosciences, Society for Industrial and Applied Mathematics
  ({SIAM}), 2015.

\bibitem{Patankar72}
S.V Patankar and D.B Spalding.
\newblock A calculation procedure for heat, mass and momentum transfer in
  three-dimensional parabolic flows.
\newblock {\em International Journal of Heat and Mass Transfer}, 15(10):1787 --
  1806, 1972.

\bibitem{Peaceman77}
D.W. Peaceman.
\newblock {\em Fundamentals of Numerical Reservoir Simulation}.
\newblock Elsevier, Amsterdam, 1977.

\bibitem{RelPerm}
PetroWiki.
\newblock Relative permeability models.
\newblock \url{http://petrowiki.org/Relative_permeability_models}.

\bibitem{RugeStueben86}
J.~W. Ruge and K.~Stueben.
\newblock {\em 4. Algebraic Multigrid}, chapter~4, pages 73--130.
\newblock 1986.

\bibitem{Russell83}
T.~F. Russell and M.F. Wheeler.
\newblock Finite element and finite difference methods for continuous flows in
  porous media.
\newblock In {\em The Mathematics of Reservoir Simulation}, chapter~2, pages
  35--106. Society for Industrial and Applied Mathematics, 1983.

\bibitem{Saad03}
Y.~Saad.
\newblock {\em Iterative Methods for Sparse Linear Systems}.
\newblock Society for Industrial and Applied Mathematics, second edition, 2003.

\bibitem{Saad86}
Y.~Saad and M.H. Schultz.
\newblock {GMRES}: A generalized minimal residual algorithm for solving
  nonsymmetric linear systems.
\newblock {\em SIAM J. Sci. Stat. Comput.}, 7(3):856--869, July 1986.

\bibitem{Stueben01}
K.~Stueben.
\newblock A review of algebraic multigrid.
\newblock {\em Journal of Computational and Applied Mathematics},
  128(1–2):281 -- 309, 2001.
\newblock Numerical Analysis 2000. Vol. VII: Partial Differential Equations.

\bibitem{Stueben07}
K.~Stueben, T.~Clees, H.~Klie, B.~Lu, and M.F. Wheeler.
\newblock Algebraic multigrid methods ({AMG}) for the efficient solution of
  fully implicit formulations in reservoir simulation.
\newblock In {\em {SPE} Reservoir Simulation Symposium}. Society of Petroleum
  Engineers ({SPE}), 2007.

\bibitem{Wallis85}
J.R. Wallis, R.P. Kendall, and L.E. Little.
\newblock Constrained residual acceleration of conjugate residual methods.
\newblock 1985.

\bibitem{Watts85}
J.W. Watts.
\newblock A compositional formulation of the pressure and saturation equations.
\newblock In {\em {SPE} Reservoir Simulation Symposium}. Society of Petroleum
  Engineers ({SPE}), SPE 12244, 1985.

\bibitem{White11}
Joshua~A. White and Ronaldo~I. Borja.
\newblock Block-preconditioned {Newton--Krylov} solvers for fully coupled flow
  and geomechanics.
\newblock {\em Computational Geosciences}, 15(4):647--659, 2011.

\bibitem{Yang06}
U.M. Yang.
\newblock Parallel algebraic multigrid methods — high performance
  preconditioners.
\newblock In {\em Numerical Solution of Partial Differential Equations on
  Parallel Computers}, volume~51 of {\em Lecture Notes in Computational Science
  and Engineering}, pages 209--236. Springer Berlin Heidelberg, 2006.

\end{thebibliography}
\end{document}